%% file: arwcourse.tex
\documentclass[a4paper,12pt,oneside]{article}

\input{macros.tex}

\begin{document}

\title{Activated Random Walks\footnote{Minicourse given at the workshop ``Activated Random Walks, DLA, and related topics'' at IMéRA-Marseille, March 2015. \textbf{Preliminary version. July 15, 2015.}}}
\author{Leonardo T. Rolla}
\date{}
\maketitle

{
\baselineskip 0pt
\tableofcontents
}

\vfill

\section{Overview}

\subsection{The Activated Random Walk reaction-diffusion model}

The Activated Random Walk model is defined as follows.

Particles sitting on the graph $\Z^d$ can be in state $A$ for \emph{active} or $S$ for \emph{passive}.
Each active particle, that is, each particle in the $A$ state, performs a continuous-time random walk with jump rate $D_A=1$ and with translation-invariant jump distribution given by a probability $p(\cdot)$ on $\Z^d$.

Several active particles can be at the same site, and they do not interact among themselves.
When a particle is alone, it may become passive, a transition denoted by $A\to S$, which occurs at a \emph{sleeping rate} $0 < \lambda \leqslant \infty$.

In other words, each particle carries two clocks, one for jumping and one for sleeping.
Once a particle is passive, it stops moving, i.e., it has jump rate $D_S=0$, and it remains passive until the instant when another particle is present at the same vertex.
At such an instant the particle which is in $S$ state flips to the $A$ state, giving the transition $A+S \to 2A$.

If the clock rings for a particle to sleep while it shares a vertex with other particles, the tentative transformation $A \to S$ is overridden by the instantaneous reaction $A+S \to 2A$, so this attempt to sleep has no effect on the system state.

\begin{framed}
Jump distribution: $x \to x+z$ with probability $p(z)$.
\\
Diffusion: jump at rates $D_A=1$ and $D_S=0$. No interaction.
\\
Reactions: $A \to S$ at rate $0 < \lambda \leqslant \infty$, $A+S \to 2A$ at rate $\infty$.
\\
Notation: $\eta_t(x)$ denotes both the number and type of particles at site $x$ at time $t$.
\end{framed}

A particle in the $S$ state stands still forever if no other particle ever visits the vertex where it is located.
At the extreme case $\lambda=+\infty$, when a particle visits an empty site, it becomes passive instantaneously.
This case is equivalent to \emph{internal diffusion-limited aggregation} with infinitely many sources.

We have described \emph{local rules} for the system to evolve.
In order to fully describe the system, we need to specify on which subset of $\Z^d$ this dynamics will occur, what are the boundary conditions, and the initial state at $t=0$.

\subsection{Infinite conservative system, fixation, phase transition}

Consider a system running on the whole graph $\Z^d$, and such that the initial configuration $\eta_0$ is i.i.d.\ Poisson with parameter $\mu$.

A phase transition in this systems arises from a conflict between a \emph{spread of activity} and a \emph{tendency for this activity to die out}, and the transition point separates an active and an absorbing phase in which the dynamics gets eventually extinct in any finite region.

\begin{framed}
We say that the system \emph{locally fixates} if $\eta_t(x)$ is eventually constant for each $x$, otherwise we say that the system \emph{stays active}.
\end{framed}

It can be shown using ergodicity that the probability that the system fixates is either $0$ or $1$.
Moreover,

\begin{framed}
{If the system fixates for given $\lambda$ and $\mu$, then it fixates for larger $\lambda$ and smaller $\mu$.}
\end{framed}

As a consequence,
\[
\Pb( \text{fixation} ) =
\begin{cases}
1 , & \mu<\mu_c, \\
0 , & \mu>\mu_c,
\end{cases}
\]
where the critical density $\mu_c=\mu_c(\lambda) \in [0,\infty]$ is non-decreasing in $\lambda \in [0,\infty]$.

\subsection{Physical motivation: ``self-organized criticality''}

In contrast with the above description, consider now the following evolution.

\textbf{Driven-dissipative system.}
For large $L$, let the system run on the finite box $V=[0,L]^d$.
New particles are added to the bulk of $V$ at constant rate.
When a particle hits the boundary of $V$ it is killed.
The reaction-diffusion dynamics is run at a much faster speed, so that the whole box is stabilized between two arrivals.
We then let $t\to\infty$ to reach stationarity, then take $L \gg 1$ to have a state on $\Z^d$.
\marginnote{Add a figure}

The relation between self-organized and ordinary criticality is understood as follows.

On the one hand, ``self-organized criticality'' appears in the parameter-free, driven-dissipative evolution described above.
In this dynamics, when the average density $\mu$ inside the box is too small, mass tends to accumulate.
When it is too large, there is intense activity and a substantial number of particles is absorbed at the boundary.
With this \emph{carefully designed} mechanism, the model is attracted to a \emph{critical state} with an average density given by $0<\mu_c<\infty$, though it was not explicitly tuned to this critical value.

On the other hand, the corresponding conservative system in infinite volume exhibits ordinary criticality in the sense that its dynamics fixate for $\mu<\mu_c$ and do not fixate for $\mu>\mu_c$, and moreover the \emph{critical exponents} of the finite-volume addition-relaxation dynamics are related to those of the conservative dynamics in infinite volume.

\subsection{Predictions}

The behavior of the ARW is expected to be the following.

\begin{figure}[b!]
\small
\centering
\includegraphics[page=2,scale=.8]{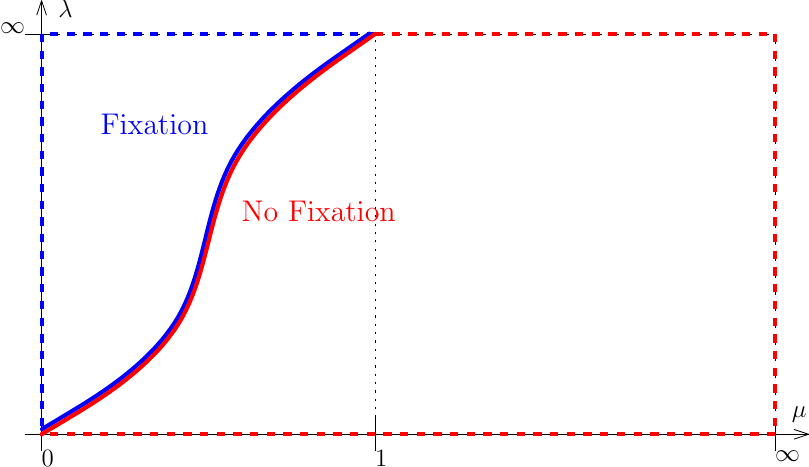}
\caption{Predictions}
\label{fig:predictions}
\end{figure}

The critical density satisfies $\mu_c \to 0$ as $\lambda \to 0$ and $\mu_c \to 1$ as $\lambda \to \infty$.
See Figure~\ref{fig:predictions}.
The value of $\mu_c$ should not depend on the particular \mbox{$\mu$-}parametrized distribution of the initial configuration (geometric, Poisson, etc.).

At $\mu=\mu_c$, the density of active particles vanishes as $t\to\infty$, but we conjecture that the system does not fixate in this case.

The asymptotic decay of density of activity as $t \gg 1$ when $\mu=\mu_c$ should obey a power law.
Also, for the stationary regime, i.e., letting $t\to\infty$ first, the density of activity should decay with a power law as $0 < \mu-\mu_c \ll 1$.
\marginnote{Include more predictions}
More predictions can be found at \cite{dickman-rolla-sidoravicius-10}.

\subsection{Difficulties}

For the ARW, we would like to describe the critical behavior, the scaling relations and critical exponents, and whether the critical density is the same as the long-time limit attained in the driven-dissipative version.
These questions are however far beyond the reach of current techniques.

The first apparent difficulty lies in the fact that this system is not attractive (i.e., its evolution does not preserve monotonicity of configurations).
This is overcome by considering a Diaconis-Fulton kind of construction, or an explicit construction in terms of a collection of random walks, rather than the classical Harris graphical construction.
These constructions allow for different kinds of arguments which have proven to be very useful.

However, another particularity of this model brings severe difficulties, which is particle conservation.
In particular, this seems to rule out any energy-entropy kind of argument.
These arguments typically go as follows.
One first identifies some structure that is intrinsic to the occurrence of events that conjecturally should not occur.
The number of possible structures is then estimated, and when it is overwhelmed by the high energy
cost needed to construct them, one can show that such events have vanishing probability.
This approach has proven successful in perturbative statistical mechanics.
However, the conservation of particles in the system gives rise to intricate long-range effects, which makes it difficult to find a suitable structure behind self-sustained activity for large periods of time.

\begin{framed}
\sout{Lack of attractiveness}.
Overcome by using other constructions rather than Harris'.

\textbf{Conservation of particles.}
Rules out any ``energy x entropy'' kind of approach.
\end{framed}

\subsection{Results}

We end this overview by recalling the existing results in the literature about the phase transition of the ARW on the infinite lattice $\Z^d$.
In the next sections we should go through the proof of all the results mentioned here.

The results will be presented in the same framework, by making the following assumptions.

\begin{framed}
We shall assume that:
(i)
the initial distribution is i.i.d.\ Poisson parametrized by $\mu$,
and
(ii)
jumps are to nearest neighbors, that is, $p(z)=0$ unless $\|z\|=1$.
\end{framed}

For each result, the assumption of i.i.d.\ initial condition and a Poisson law for its marginals can be more or less relevant depending on the proof.
We cannot make any other assumptions about the jump distribution, since some of the existing results have been proved only for biased jumps, whereas some others only in the complementary case.

Each result can be generalized in different directions.
The reader interested in minimal hypotheses is referred to the original articles.

\subsubsection*{One-dimensional directed walks}

For the particular situation when particles only jump to their right, all the above predictions about the phase transition curve can be proved.
More precisely,
\[
\mu_c = \frac{\lambda}{1+\lambda},
\]
and there is not fixation at $\mu=\mu_c$.
This appeared in \cite{cabezas-rolla-sidoravicius-14} although it was known before.
Yet, much remains to be understood about its critical behavior.
In particular, the scaling limit of the flow process has been studied in~\cite{cabezas-rolla-sidoravicius-14} for the extreme case $\lambda=\infty$, but remains open elsewhere on the critical curve.
\marginnote{Describe better the scaling limits and conjectures}
\begin{figure}[thb]
\centering
\includegraphics[page=3,scale=.8]{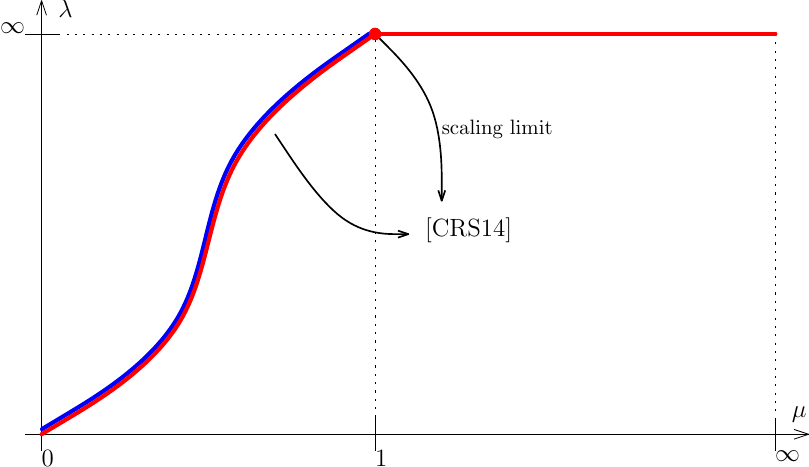}
\caption{Results for $d=1$, totally asymmetric jumps}
\end{figure}

\subsubsection*{One-dimensional systems}

In studying a different model, arguments found in~\cite{kesten-sidoravicius-06} imply that for high enough $\mu$ the ARW with $\lambda=\infty$ model stays active.
In terms of the phase diagram, $\mu_c<\infty$ for any $\lambda \leqslant \infty$.
It was shown in~\cite{rolla-sidoravicius-12} that there is no fixation at $\mu=1$, so in particular $\mu_c \leqslant 1$ for any $\lambda \leqslant \infty$.

Finally, in~\cite{taggi-14} it was shown that when the jump distributions is biased, $\mu_c<1$ for any $\lambda<\infty$, and $\mu_c\to 0$ as $\lambda \to 0$.
It should be noted that the lower bound obtained for $\mu_c$ tends to $1$ if the bias is small, for any $\lambda>0$, and that the case of unbiased jumps remains open.

The first result on fixation in this context appeared in~\cite{rolla-sidoravicius-12}, showing that $$\mu_c \geqslant \frac{\lambda}{1+\lambda}.$$
In terms of the phase diagram, $\mu_c>0$ for all $\lambda$, and $\mu_c\to 1$ as $\lambda \to \infty$.

\begin{figure*}[htb]
\small
\centering
\includegraphics[page=4,scale=.8]{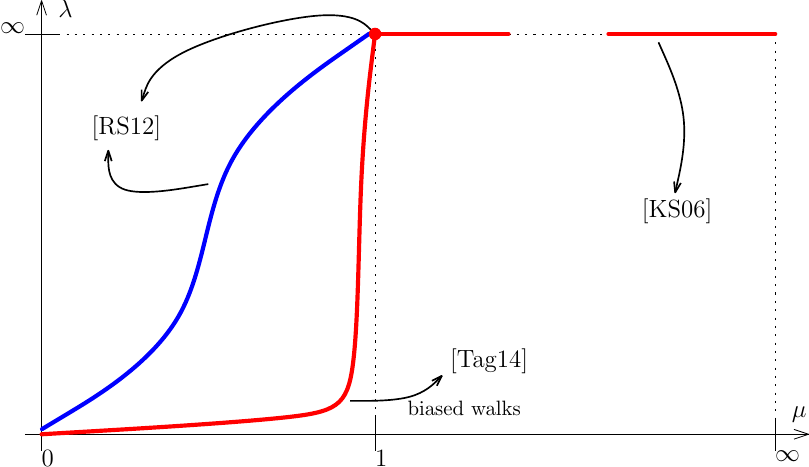}
\caption{Results for $d=1$, general jumps}
\end{figure*}

\subsubsection*{Two and more dimensions}

In~\cite{kesten-sidoravicius-06}, the fact that $\mu_c < \infty$ was proved for any dimension.
It was then shown that $\mu_c \leqslant 1$ by~\cite{shellef-10} and~\cite{amir-gurelgurevich-10}, using different methods.
This was extended in~\cite{cabezas-rolla-sidoravicius-14}, who showed that there is no fixation at $\mu=1$.

In~\cite{taggi-14}, it was shown that $\mu_c<1$ for sufficiently biased jump distributions and some $\lambda>0$.
Unlike the one-dimensional case, the proof gives neither $\mu_c<1$ for all $\lambda<\infty$, nor $\mu_c \to 0$ as $\lambda \to 0$.
The case of small or zero bias remains open.

The first result on fixation for arbitrary dimensions appeared in~\cite{shellef-10}, who shows that $\mu_c>0$ when $\lambda=\infty$.
This was extended in~\cite{cabezas-rolla-sidoravicius-14} who show that $\mu_c = 1$ when $\lambda=\infty$.

Finally, in~\cite{sidoravicius-teixeira-14}, it was shown that when the jump distribution is unbiased, $\mu_c>0$ for any $\lambda>0$.
The methods of~\cite{sidoravicius-teixeira-14} do not allow to prove that $\mu_c \to 1$ as $\lambda \to \infty$, though.
The case of biased jumps and $\lambda<\infty$ remains widely open.

\begin{figure*}[hbt]
\small
\centering
\includegraphics[page=5,scale=.8]{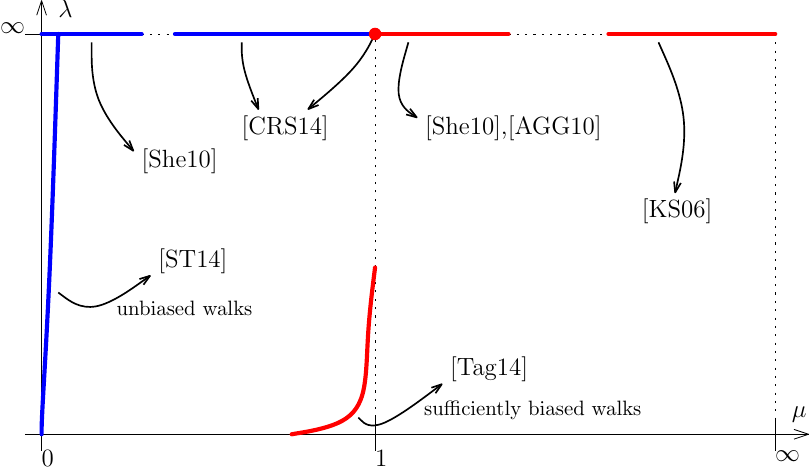}
\caption{Results for $d \geqslant 2$}
\label{fig:highdpredic}
\end{figure*}

\section{Definitions and Diaconis-Fulton construction}

\subsection{Notation}

Let $\N_0 = \{0,1,2,\dots\}$ and $\N_{\varrho} = \N_0 \cup \{\varrho\}$.
The state of the ARW at time $t\geqslant 0$ is given by $\eta_t \in \Sigma = (\N_{\varrho})^{\Z^d}$.
In this setting, $\eta_t(x)=n$ means that, at time $t$, site $x$ contains $n$ active particles if $n\in\N_0$, or one passive particle if $n=\varrho$.

We turn $\N_{\varrho}$ into an ordered set by letting $0<\varrho<1<2<\cdots$.
We also let $|\varrho|=1$, so $|\eta_t(x)|$ counts the number of particles regardless of their state.
We define $\llbracket n \rrbracket$ to be $n$ if $n \geqslant 1$, and $0$ if $n=0 \text{ or } \varrho$, so $\llbracket \eta_t(x) \rrbracket$ counts the number of active particles.

To add a particle to a site we define $\varrho+1=2$, which represents the $A+S\to2A$ transition.
We also define $1 \cdot \varrho = \varrho$ and $n \cdot \varrho = n$ for $n \geqslant 2$, which represent the transitions $A \to S$ and $2A \to A+S \to 2A$, respectively.

It will be convenient to consider as \emph{acceptable} the operations $\varrho-1=0$ and $\varrho\cdot\varrho=\varrho$, although they do not appear in the dynamics.
The operations $0\cdot\varrho$ and $0-1$, on the other hand, are not acceptable.

The process evolves as follows.
For each site $x$, a Poisson clock rings at rate $(1+\lambda) \llbracket \eta_t(x) \rrbracket$.
When this clock rings, the system goes through the transition $\eta\to \tau_{x\varrho}\eta$ with probability $\frac{\lambda}{1+\lambda}$, otherwise $\eta\to \tau_{xy}\eta$ with probability $p(y-x)\frac{1}{1+\lambda}$.
The transitions are given by
\[
  \tau_{xy}\eta(z) =
  \begin{cases}
    \eta(x)-1, & z=x \\
    \eta(y)+1, & z=y \\
    \eta(z),   & \mbox{otherwise,}
  \end{cases}
\qquad
\text{and}
\qquad
  \tau_{x\varrho}\eta(z) =
  \begin{cases}
    \varrho\cdot\eta(x), & z=x \\
    \eta(z),   & \mbox{otherwise}
.
  \end{cases}
\]

We assume that $\eta_0(x)\in\N_0$ for all $x$ a.s., and use $\Pb^\nu$ denote the law of $(\eta_t)_{t\geqslant0}$, where $\nu$ denotes the distribution of $\eta_0$.

\subsection{Diaconis-Fulton construction}

The Diaconis-Fulton representation enables us to exploit the combinatorial nature of fixation.
Due to particle exchangeability, this representation extracts precisely the part of the randomness that is relevant for the phase transition, focusing on the total number of jumps and leaving aside the order in which they take place.
It is suitable for studying path traces, total occupation times, and final particle positions, but precludes the analysis of quantities for which the order of the jumps does matter, such as correlation functions or local shape properties.
\marginnote{Explain what the DF construction is and add a picture}

We summarize the main points of this section for later reference.

\begin{framed}
--
\emph{Toppling $x$} is an operator that reduces the state at $x$ and increases elsewhere.
\\
--
\emph{Legal topplings} are those actually performed by the system.
\\
--
The order of topplings is irrelevant and we are free to choose it.
\\
--
\emph{Acceptable toppling} is an artificial operation that extends the definition of legal topplings, and provides \emph{upper bounds} for activity.
\\
--
On the other hand, legal topplings provide \emph{lower bounds} for activity.
\\
--
This construction determines fixation or non-fixation for the continuous-time model.
\\
--
There is a zero-one law for fixation.
\end{framed}

In this section, $\eta$ denotes a configuration, and we do not deal with a continuous-time evolution anymore.
We say that site $x$ is \emph{unstable for the configuration $\eta$} if $\eta(x) \geqslant 1$.
Unstable sites can \emph{topple}, by applying $\tau_{xy}$ or $\tau_{x\varrho}$ to $\eta$.
Toppling an unstable site is \emph{legal}.
If $\eta(x) = \varrho$, we still say that toppling $x$ is \emph{acceptable}.
A legal toppling is also acceptable.

Consider a field of instructions $\mathcal{I}=(\tau^{x,j})_{x\in\Z^d, j\in\N}$.
Later on we will choose $\II$ random, but for now $\mathcal{I}$ denotes a field that is fixed, and $\eta$ denotes any configuration.

Let $\h=\big(h(x);x\in\Z^d\big)$ count the number of topplings at each site, usually started from $h\equiv 0$.
The toppling operation at $x$ is defined by
\[
\Phi_x(\eta,\h)= \big(\tau^{x,h(x)+1}\eta, \h + \delta_x \big)
,
\]
and
$\Phi_x\eta$ is a short for $\Phi_x(\eta,0)$.
Again, $\Phi_x$ is \emph{legal for $\eta$} if $\eta(x) \geqslant 1$, and \emph{acceptable for $\eta$} if $\eta(x) \geqslant \varrho$.
\marginnote{Add a figure}

\subsubsection*{Sequences of topplings and local properties}

Let $\alpha=(x_1,\dots,x_k)$ denote a finite sequence of sites.
We define $\Phi_\alpha = \Phi_{x_k}\Phi_{x_{k-1}}\cdots\Phi_{x_1}$.
We say that $\Phi_\alpha$ is \emph{legal} (respectively \emph{acceptable}) for $\eta$ if $\Phi_{x_l}$ is legal (respectively acceptable) for $\Phi_{(x_{1},\dots,x_{l-1})}\eta$ for each $l=1, \dots, k$.
In this case we say that $\alpha$ is a \emph{legal sequence} (respectively \emph{acceptable sequence}) of topplings for $\eta$.

Let $m_\alpha=\big(m_\alpha(x);x\in\Z^d\big)$ be given by
$m_\alpha(x)=\sum_l{\I}_{x_l=x}$, the number of times the site $x$ appears in
$\alpha$.
We write $m_\alpha \geqslant m_\beta$ if $m_\alpha(x) \geqslant m_\beta(x)\ \forall\ x$, and $\tilde{\eta} \geqslant \eta$ if $\tilde{\eta}(x) \geqslant \eta(x)\ \forall\ x$.
We also write $(\tilde{\eta},\tilde h) \geqslant (\eta,h)$ if $\tilde{\eta} \geqslant \eta$ and $\tilde{h} = h$.

Let $x$ be a site in $\Z^d$ and $\eta,\tilde\eta$ be configurations.

\begin{framed}
\textbf{1) Local Abelianness.}
If $\alpha$ and $\beta$ are acceptable sequences of topplings for the configuration $\eta$, such that $m_\alpha=m_\beta$, then $\Phi_\alpha \eta = \Phi_\beta \eta$.

\textbf{2) Mass comes from outside.}
If $\alpha$ and $\beta$ are acceptable sequences of topplings for $\eta$ such that $m_\alpha(x) \leqslant m_\beta(x)$ and $m_\alpha(z) \geqslant m_\beta(z)$ for $z\ne x$, then $\Phi_\alpha\eta(x) \geqslant \Phi_\beta\eta(x)$.

\textbf{3) Monotonicity of stability.}
If site $x$ is unstable for the configuration $\eta$, and if $\tilde\eta(x)\geqslant\eta(x)$, then $x$ is also unstable for the configuration $\tilde\eta$.

\textbf{4) Monotonicity of topplings.}
If $\tilde\eta\geqslant\eta$ and $\Phi_x$ is legal for $\eta$, then $\Phi_x$ is legal for $\tilde\eta$ and $\Phi_x \tilde\eta \geqslant \Phi_x \eta$.
\end{framed}

\begin{proof}
The two last properties are immediate.
Suppose that $m_\alpha=m_\beta$.

For convenience, define the operators $n \oplus = n+1$ on $\N_\varrho$, we well as $n \ominus = n-1$ and $n\odot = n \cdot \varrho$ on $\N_\varrho \setminus \{0\}$.
With this notation, whenever $n \ominus$ is acceptable (i.e., $n\ne 0$) we have $n \ominus \oplus = n \oplus \ominus$.
Analogously, whenever $n \odot$ is acceptable (i.e., $n\ne 0$) we have $n \odot \oplus = n \oplus \odot$.
Therefore, within any acceptable sequence of operations, replacing $\ominus \oplus$ and $\odot \oplus$ by $\oplus \ominus$ and $\oplus \odot$ yields an acceptable sequence with the same final outcome.

Notice that $\Phi_\alpha \eta (x)$ is given by $\eta(x)$ followed by a number of $\oplus$, $\ominus$ and $\odot$'s.
The number of times each operator appears is determined by $\II$ and $m_\alpha$ and is thus the same for $\Phi_\beta \eta (x)$.
Their actual order depends on the sequence, but the internal order of the $\ominus$ and $\odot$ operators is determined by $(\tau^{x,j})_j$ and is the same for both $\Phi_\alpha \eta (x)$ and $\Phi_\beta \eta (x)$.
As a consequence, we can apply the above identities to move the $\oplus$'s to the left, yielding then identical sequences for $\Phi_\alpha \eta (x)$ and $\Phi_\beta \eta (x)$, which proves the first property.

The second property follows from a similar argument.
Suppose $m_\alpha(x) \leqslant m_\beta(x)$ and $m_\alpha(z) \geqslant m_\beta(z)$ for $z\ne x$.
Again $\Phi_\alpha \eta (x)$ is given by $\eta(x)$ followed by a number of $\oplus$, $\ominus$ and $\odot$'s.
The number of times that operator $\oplus$ appears depends on $m_\alpha(z), z\ne x$, and is bigger than for $\Phi_\beta \eta (x)$, and the number of times that operators $\ominus$ and $\odot$ appear depend on $m_\alpha(x)$, and is smaller than for $\Phi_\beta \eta (x)$.
Pushing the $\oplus$'s to the left we get that $\Phi_\alpha \eta(x)$ is written in the same way as $\Phi_\beta \eta(x)$, perhaps with a few extra $\oplus$'s in the beginning, and a few missing $\ominus$ and $\odot$'s in the end, so $\Phi_\alpha\eta(x) \geqslant \Phi_\beta\eta(x)$.
\end{proof}

\subsubsection*{Global properties}

What follows is valid for any model satisfying the above four properties.

Let $V$ denote a finite subset of $\Z^d$.
We say that $\eta$ is \emph{stable in $V$} if every $x\in V$ is stable for $\eta$.
We say that $\alpha$ is \emph{contained in $V$}, and write $\alpha \subseteq V$, if every $x$ appearing in $\alpha$ is an element of $V$.
We say that $\alpha$ \emph{stabilizes $\eta$ in $V$} if $\alpha$ is acceptable for $\eta$ and $\Phi_\alpha\eta$ is stable in $V$.

\begin{framed}
\textbf{Least Action Principle.}
If $\alpha$ is an acceptable sequence of topplings that stabilizes $\eta$ in $V$, and $\beta \subseteq V$ is a sequence of topplings that is legal for $\eta$, then $m_\beta \leqslant m_\alpha$.
\end{framed}
\begin{proof}
Let $\beta \subseteq V$ be legal and $m_\alpha \ngeqslant m_\beta$.
Write $\beta=(x_1,\dots,x_k)$ and $\beta^{(j)}=(x_1,\dots,x_j)$ for $j\leqslant k$.
Let $\ell=\max\{ j : m_{\beta^{(j)}} \leqslant m_{\alpha} \}<k$ and $y=x_{\ell+1}\in V$.
Since $\beta$ is legal, $y$ is unstable in $\Phi_{\beta^{(\ell)}}\eta$.
But $m_{\beta^{(\ell)}}\leqslant m_{\alpha}$ and $m_{\beta^{(\ell)}}(y)= m_{\alpha}(y)$.
By the Properties~2 and~3, $y$ is unstable for $\Phi_{\alpha}\eta$ and therefore $\alpha$ does not stabilize $\eta$ in $V$.
\end{proof}

For each finite set $V$ and configuration $\eta$, we define
\[
m_{V,\eta} = \sup_{\beta\subseteq V\text{ legal}} m_\beta.
\]
The Least Action Principle says that
\[
m_\beta \leqslant m_{V,\eta} \leqslant m_\alpha
\]
for any legal sequence $\beta$ contained in $V$ and any acceptable sequence $\alpha$ stabilizing $\eta$ in $V$.
This provides a very good source of lower and upper bounds for $m_{V,\eta}$.

\begin{framed}
\textbf{Global Abelianness.}
If $\alpha$ and $\beta$ are both legal toppling sequences for $\eta$ that are contained in $V$ and stabilize $\eta$ in $V$, then $m_\alpha=m_\beta=m_{V,\eta}$.
In particular, $\Phi_\alpha\eta=\Phi_\beta\eta$.
\end{framed}
\begin{proof}
Apply the Least Action Principle in two directions: $m_\beta \leqslant m_{V,\eta} \leqslant m_\alpha \leqslant m_\beta$.
\end{proof}

We say that \emph{$\eta$ is stabilizable in $V$ } if there exists an acceptable sequence $\alpha$ that stabilizes $\eta$ in $V$.
Notice that this provides finite upper bounds for legal sequences $\beta\subseteq V$, which in turn implies the existence of a sequence contained in $V$ that is both legal and stabilizing.
Indeed, if one tries to perform legal topplings in $V$ indefinitely, on the one hand one has to eventually stop since there is a finite upper bound, and on the other hand this can only be stopped if there are no more unstable sites.

\begin{framed}
\textbf{Monotonicity.}
If $V\subseteq \tilde V$ and $\eta \leqslant \tilde\eta$, then $m_{V,\eta}\leqslant m_{\tilde V,\tilde\eta}$.
\end{framed}
\begin{proof}
Let $\beta \subseteq V$ be legal for $\eta$.
By successively applying Properties~3 and~4, $\beta$ is also legal for $\tilde\eta$.
Since $\beta \subseteq \tilde V$ in this case, the inequality follows from the definition of $m_{V,\eta}$.
\end{proof}

By monotonicity, the limit
\[
m_\eta = \lim_{V \uparrow\Z^d} m_{V,\eta}
\]
exists and does not depend on the particular sequence $V \uparrow\Z^d$ (the limit is given by the supremum over finite $V$).
A configuration $\eta$ is said to be \emph{stabilizable} if $m_\eta(x)<\infty$ for every $x\in\Z^d$.

\subsubsection*{Fixation for the stochastic dynamics}

Recall that $\Pb^\nu$ denotes the law of the process $(\eta_t)_{t\geqslant 0}$ with values on $\N_\varrho^{\Z^d}$.
Assume that distribution $\nu$ of $\eta_0$ on $\N_0^{\Z^d}$ is ergodic and has finite density $\nu(\eta(\o))<\infty$.
We further assume that the support of $p(\cdot)$ generates $\Z^d$ and not a sublattice.

To avoid extra notation, we  define the field $\II$ on the same probability space $\Pb^\nu$, independently of $\eta_0$, and distributed as follows.
For each $x\in\Z^d$ and $j\in\N$, choose $\tau^{x,j}$ as $\tau_{xy}$ with probability $\frac{p(y-x)}{1+\lambda}$ or $\tau_{x\varrho}$ with probability $\frac{\lambda}{1+\lambda}$, independently over $x$ and $j$.

\begin{framed}
\textbf{Stochastic Dynamics.}
$\Pb^\nu(\text{fixation of }(\eta_t)_{t\geqslant 0})=\Pb^\nu (m_{\eta_0}(\o)<\infty) = 0 \text{ or } 1$.
\end{framed}

For the proof, the reader is referred to~\cite[v2 on arXiv]{rolla-sidoravicius-12}.

\section{One-dimensional counting arguments}

\begin{framed}
\cite{cabezas-rolla-sidoravicius-14} ($d=1$. Directed walks).
$\mu_c=\frac{\lambda}{1+\lambda}$.
No fixation for $\mu=\mu_c$.
\end{framed}

\begin{proof}
Fixation is equivalent to $m_{\eta_0}(\o)<\infty$ almost surely.
This is equivalent to $m_{V_L,{\eta_0}}(\o)$ being tight as $L\to\infty$, where $V_L=\{-L,\dots,0\}$.
The latter is equivalent to tightness of the number of particles which jump into $\o$ when $V_L$ is stabilized with legal topplings.
See Figure~\ref{fighowmany}.
\begin{figure}[!htb]
\centering
\includegraphics[page=1,scale=1]{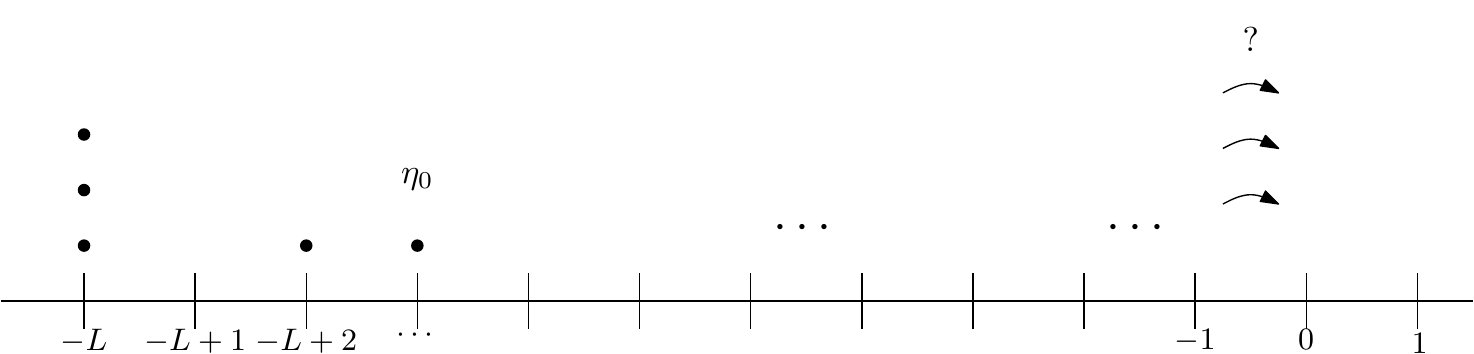}
\caption{The set $V_L$, part of the configuration $\eta$, and arrow indicating the unknown number of particles which jump into $\o$ when stabilizing $[-L,0]$}
\label{fighowmany}
\end{figure}

Let us stabilize $V_L$ by exhausting each site, from left to right.
More precisely, topple site $x=-L$ until each of the $\eta_0(-L)$ particles either moves to $x=-L+1$ or sleeps, and let $Y_0$ be indicator of the event that the last particle remained sleeping on $x=-L$.
Conditioned on $\eta_0(-L)$, the distribution of $Y_0$ is Bernoulli with parameter $\frac{\lambda}{1+\lambda}$ (in case $\eta_0(-L)=0$, sample $Y_0$ independently of anything else).
The number of particles which jump from $x=-L$ to $x=-L+1$ is given by $N_1:=[\eta_0(-L)-Y_0]^+$.
See Figure~\ref{figcounting}.
\begin{figure}[!htb]
\centering
\includegraphics[page=2,scale=1]{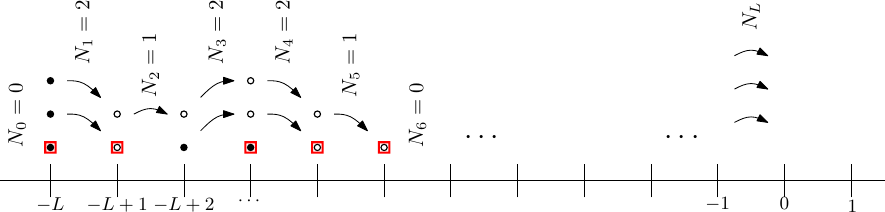}
\caption{Stabilizing $[-L,0]$ from left to right. Disks represent particles initially present at each site, and circles represent particles that arrive from the left. Red boxes indicate particles which became passive when they were left alone.}
\label{figcounting}
\end{figure}

Note that, after stabilizing $x=-L$, there are $N_1+\eta_0(-L+1)$ particles at $x=-L+1$.
Let the site $x=-L+1$ topple until it is stable, and denote by $Y_1$ the indicator of the event that the last particle remained passive on $x=-L+1$.
The number of particles which jump from $x=-L+1$ to $x=-L+2$ is given by $N_2:=[N_1+\eta_0(-L+1)-Y_1]^+$.
By iterating this procedure, the number $N_{i+1}$ of particles which jump from $x=-L+i$ into $x=-L+i+1$ after stabilizing $x=-L,-L+1,\dots,-L+i$ is given by $$N_{i+1}=[N_{i}+\eta_0(-L+i)-Y_i]^+,$$ where $N_0=0$.
The number of particles which jump into $\o$ while stabilizing $V_L$ equals $N_L$.

Note that the sequence $(N_i)_{i=0,1,\dots,L}$ is a random walk on $\N_0$, with independent increments distributed as ${\eta_0}(x)-Y$, reflected at $0$.
So the relevant quantity is
$$
\E[{\eta_0}(-L+k)-Y_k]=\mu - \frac{\lambda}{1+\lambda}
.
$$
If $\mu < \frac{\lambda}{1+\lambda}$, the walk is positive recurrent, which implies tightness of $N_{i}$, and thus fixation.
If $\mu > \frac{\lambda}{1+\lambda}$, the walk is transient, and $\Pb[N_{L} \geq \frac{1}{2}(\mu - \frac{\lambda}{1+\lambda})L]$ is large for large $L$, so there is no fixation.
If $\mu = \frac{\lambda}{1+\lambda}$, the walk is null-recurrent, so there is no tightness and thus no fixation.
\end{proof}

\pagebreak[2]
\begin{framed}
($d=1$).
No fixation when $\mu=1$.
\end{framed}

\begin{proof}
Let $\mu=1$.
By the CLT, the probability that $\eta_0$ contains at least $L+2 \sqrt{L}$ particles on $V_L=[0,L]$ is at least $2 \delta>0$, uniformly on $L$.
On this event, at least $2\sqrt{L}$ particles will visit either $x=\o$ or $x=L$ when we stabilize $[0,L]$.
Therefore, $$2\delta \leqslant \Pb\big(m_{V_L,{\eta_0}}(\o) \geqslant \sqrt{L}\big) + \Pb\big(m_{V_L,{\eta_0}}(L) \geqslant \sqrt{L}\big) \leqslant 2 \Pb(m_{\eta_0}(\o) \geqslant \sqrt{L}).$$
Since this is true for any $L$, there is no fixation.
\end{proof}

\begin{framed}
\cite{taggi-14} ($d=1$, Biased walks).
For any $\lambda<\infty$, $\mu_c<1$. Moreover, $\lim_{\lambda \to 0}\mu_c(\lambda) = 0$.
\end{framed}

\begin{proof}
Let $F(\lambda)$ be the probability that a walk which jumps at rate 1 and sleeps at rate $\lambda$, but only sleeps at $(-\infty,\o]$, is able to reach $+\infty$ without sleeping.
A more explicit definition is $F(\lambda)= \E [(\frac{1}{1+\lambda})^T ]$, where $T$ counts how many steps a discrete-time random walk starting at $x=\o$ is found on $(-\infty,0]$.
We will show that $\mu_c \leqslant 1-F(\lambda)$.

\begin{figure}[!b]
\centering
\includegraphics[page=2,scale=.8]{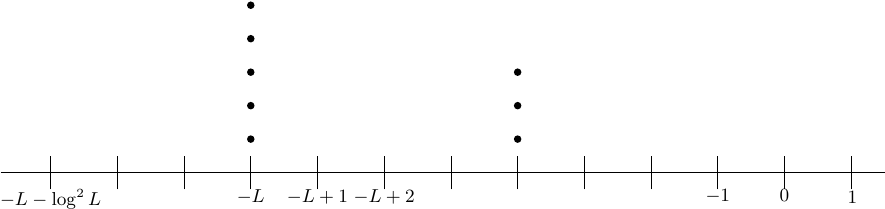}

\vspace*{5mm}
\includegraphics[page=3,scale=.8]{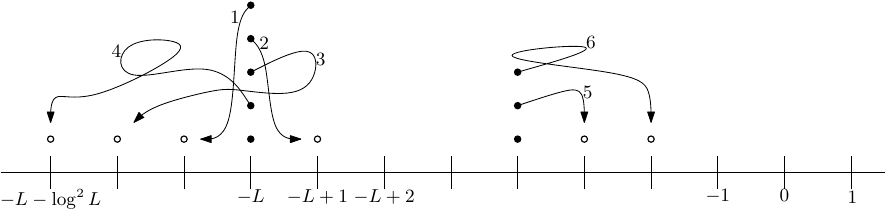}

\vspace*{5mm}
\includegraphics[page=4,scale=.8]{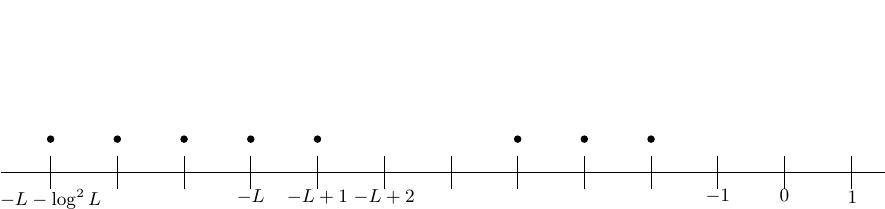}

\vspace*{5mm}
\includegraphics[page=5,scale=.8]{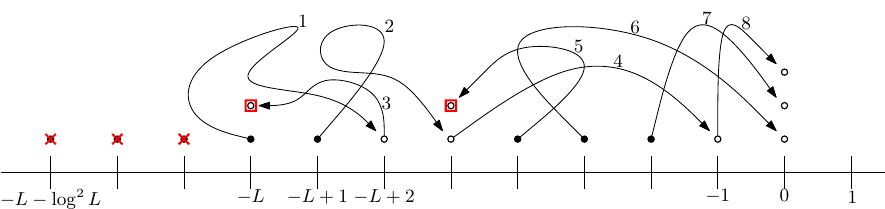}

\vspace*{5mm}
\includegraphics[page=6,scale=.8]{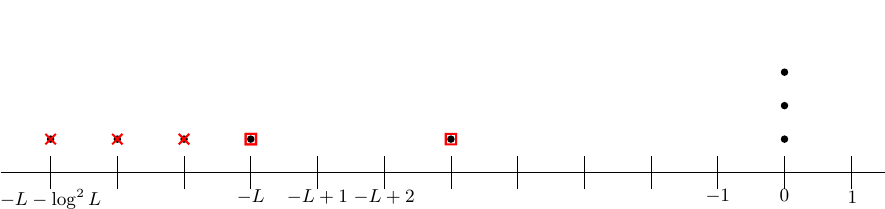}
\caption{The two stages in the toppling procedure. After the second stage, a number of particles is found at $\o$.}
\label{fig:taggi1d}
\end{figure}

Suppose $\mu>1-F(\lambda)$.
We will perform a sequence of legal topplings on $[-2L,0]$ as follows.

In the first stage, we move each particle initially found on $[-L,-1]$ until it is alone, in order to have a configuration on $\{0,1\}^L$.
If particles reach $\{0\}$ they are stopped, but on $[-2L,-L]$ we still move them.
With high probability, none of the biased random walks will go further to the left than $\log^2 L$.
Therefore, after this step at most $\log^2 L$ particles will be outside $[-L,0]$.

We now move to the second stage.
Let $N_0$ denote the number of particles on $A_0=[-L,0]$ after this first step. Notice that
\[
\text{At this point, all sites in this set, except site $\o$, have either 0 or 1 particle}
\eqno{(*)}
\]
Let $x=-L$ and $A_1=[-L+1,0]$.
If $x$ has no particles, then at $A_1$ there are $N_1=N_0$ particles, and $A_1$ satisfies (*). If $x$ has one particle, move this particle until it either sleeps on $[-\infty,x]$, or finds an empty site on $[x+1,-1]$, or else it reaches $\{0\}$.
At the end of this step, the number of particles present at $A_1$ will be either $N_0$ or $N_0-1$.
The latter happens with probability at most $1-F(\lambda)$.
In any case the configuration on $A_1$ satisfies~(*).

Move to $x=-L+1$, define $A_2=[-L+2,0]$, again move the particle found at $x$ (if any) until it reaches an empty site.
Again the number $N_2$ of particles in $A_2$ after this step is $N_1-1$ with probability at most $1-F(\lambda)$, otherwise is equal to $N_1$, and $A_2$ satisfies~(*).
Then move to $x=-L+2$, define $A_3$, $N_3$, and so on.

Finally, when $x=-1$ we have the set $A_L=\{0\}$ containing $N_L$ particles.
By applying the LLN, both to the number of particles initially found on $[-L,0]$, and to the number of particles lost on the second stage, we have that $N_L > (\mu-\epsilon) L - \log^2 L - [1-F(\lambda)+\epsilon] L \gg 1$ with high probability, finishing the proof.
\end{proof}
Notice that the upper bound for $\mu_c$ degenerates as the jumps become less biased.

\section{Exploring the instructions in advance}

\begin{framed}
\cite{rolla-sidoravicius-12} ($d=1$).
For any $\lambda>0$, $\mu_c \geqslant \frac{\lambda}{1+\lambda}$.
\end{framed}

\begin{proof}
[Sketch of the proof]
The proof uses an algorithm that tries to stabilize all the particles initially present in $\eta_0$, following the instructions in $\II$, with the aid of acceptable topplings.

After describing the algorithm, we will show that, whenever it is successful, it implies that $m_{\eta_0}(\o)=0$.
We finally show that the algorithm is successful with positive probability if $\mu < \frac{\lambda}{1+\lambda}$.
By the $0$-$1$ law, this implies almost sure fixation of the ARW model.

I.\ Description of the algorithm.

The algorithm consists in applying a trapping procedure to each particle.
This procedure explores $\II$ until it identifies a suitable trap (given by a carefully chosen sleep instruction) for the particle.
To do that, it follows the path that the particle would perform if we always toppled the site it occupies, and stops when the trap has been chosen.
In the absence of a suitable trap, the algorithm fails.

Remark that some of the explored instructions are actually not going to be used by the corresponding particle by the time it settles at the trap, leaving some \emph{corrupted sites} that could interfere with the statistics of the subsequent steps.
\marginnote{Make the pictures more descriptive}

If there is a particle at $\o$, we declare the procedure unsuccessful and stop.
Otherwise, label the initial positions of the particles on $\Z$ by $\cdots \leqslant x_{-3} \leqslant x_{-2} \leqslant x_{-1} < 0 < x_1 \leqslant x_2 \leqslant x_3 \leqslant \cdots$.
Let $a_0=0$.

We now describe the trapping procedure for particle $x_k$, $k \geqslant 1$.
Suppose that the first $k-1$ traps have been successfully set up at positions $0<a_1<a_2<\cdots<a_{k-2}<a_{k-1}<x_{k-1}$.
The set $[ 0,a_{k-1} ]$ contains all traps and corrupted sites
found so far.

The settling procedure starts with an \emph{exploration}.
Starting at $x_k$, examine and follow the instructions in $\II$ one by one
(whenever a sleep instruction is found, the next instruction
at the same site is to be examined).
Follow this exploration until reaching $a_{k-1}$.

Next we \emph{set up the trap}.
During the $k$-th exploration, we are sure to visit every $x\in
B_k=[ a_{k-1}+1,x_k-1 ]$.
Moreover, the last instruction explored at each $x\in B_k$ is a jump to the
left, see Figures~\ref{fig1arw}~and~\ref{fig2arw}.
For each $x\in B_k$, the second last instruction may be a sleep instruction.
If this is not the case for any $x\in B_k$, we declare the procedure unsuccessful and stop.
Otherwise, let $a_k$ be the leftmost site at which the second last instruction explored was a sleep instruction, and call this second last instruction the \emph{$k$-th trap}.

Since the trap is a sleep instruction found immediately before the last
instruction, which is a jump to the left, we are sure that the exploration path
has not gone to the right of $a_k$ after finding the trap, so \emph{all the
corrupted sites will be in $[ a_{k-1}+1,a_k]$}, see
Figure~\ref{fig2arw}.
This procedure can be carried on indefinitely, as long as \emph{all the steps
are successful}.

\begin{figure}[tb!]
 \centering
 \includegraphics{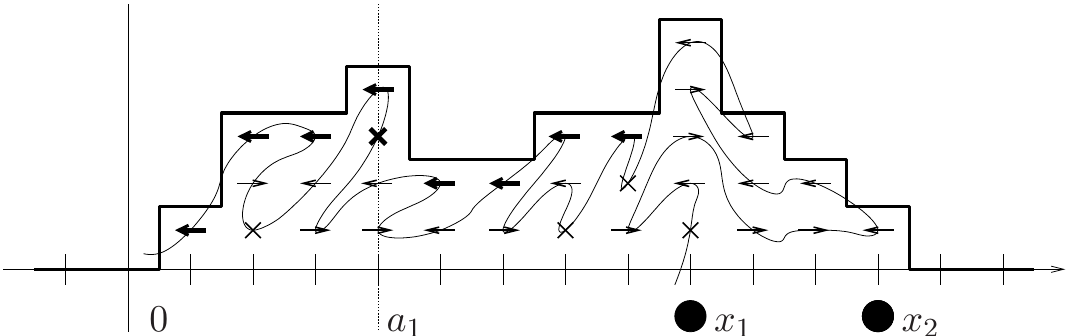}
 \caption{
 First exploration path.
 It starts at position $x_1$ of the first particle and stops when it reaches the origin.
 The horizontal axis represents the lattice, and above each site $x$ there is a sequence of instructions $(\tau^{x,j})_j$.
 The bold arrows indicate the last jump found at each site $x \in [
1,x_1-1 ]$, and the bold cross indicates a sleep instruction found just
before the last jump, this being the leftmost such cross, whose
location defines $a_1$.
 }
 \label{fig1arw}
\end{figure}

\begin{figure}[b!]
 \centering
 \includegraphics{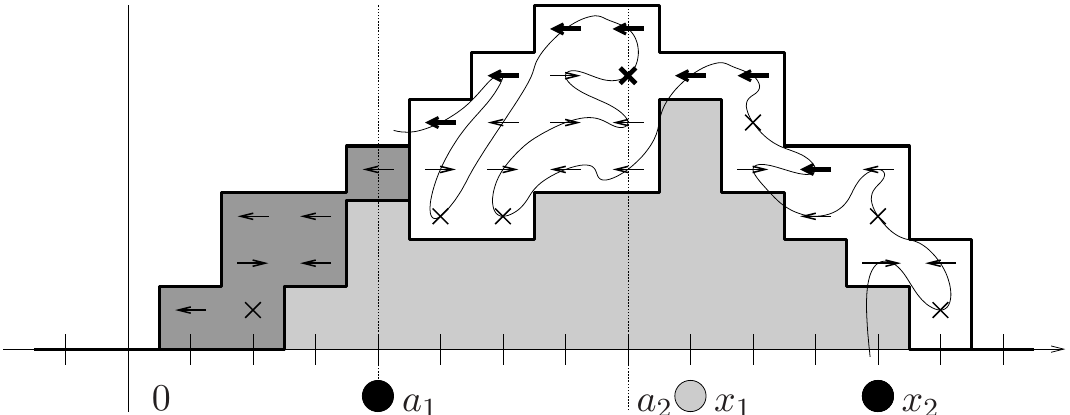}
 \caption{
 Second exploration path.
 It starts at position $x_2$ of the second particle and stops when it reaches $a_1$.
 The regions in gray indicate the instructions already examined by the first
explorer.
 The dark gray contains instructions examined but not used, whose locations
determine the set of corrupted sites.
 }
 \label{fig2arw}
\end{figure}

\medskip

II.\ If the algorithm is successful, then $m_{\eta_0}(\o)=0$.

We will show that, following the instructions of $\II$, $\eta_0$ is stabilized in $V_n=[x_{-n},x_n]$ with finitely many acceptable topplings, without toppling $\o$.

Let us first stabilize the particle that starts at $x_1$.
To that end, we successively topple the sites found by the first explorer, until it reaches the trap at $a_1$.
At this moment the particle will become passive, and the site $a_1$ will be stable.
Notice that this is acceptable since the particle is following the same path that the explorer did, even if it will sometimes be passive.

Notice that, after the last visit to $a_1$, the explorer did not go further to the right, so when settling the first particle we use all the instructions examined so far, except some lying in $[a_0+1,a_1]$.
Therefore, the same procedure can be applied to the second particle, as it will find the same instructions that determined the second exploration path.

Notice also that the first particle does not visit $\o$, and the second particle neither visits $\o$ nor $a_1$, thus it is settled without activating the first particle.
Following the same procedure, the $k$-th particle is settled at $a_k$, without ever visiting $\{0,a_1,a_2,\dots,a_{k-1}\}$, for all $k=1,\dots,n$.
After settling the $n$ first particles in $\Z_+$, we perform the analogous procedure for the first $n$ particles in $\Z_-$.

This means that $\eta_0$ can be stabilized in $V_n$ with finitely many acceptable topplings, not necessarily in $V_n$, and never toppling the origin.
By the Least Action Principle, $m_{V_n,\eta}(\o)=0$.
Since it holds for all $n\in\N$ and $V_n\uparrow\Z$ as $n\to\infty$, this gives
$m_{\eta}(\o)=0$.

\medskip

III.\ The algorithm is successful with positive probability.

For each site $x \in B_1$, the probability of finding a sleep instruction just before its last jump equals $\frac{\lambda}{1+\lambda}$, and this happens independently of the path and independently for each site.
Thus, $a_1-a_0$ is a random variable distributed as a geometric with parameter $\frac{\lambda}{1+\lambda}$, truncated at $x_1-a_0$.

Since no corrupted sites were left outside $[ a_0+1,a_1 ]$ in the first exploration, the interdistance $a_2-a_{1}$ is independent of $a_1$.
Moreover, its distribution is also geometric with parameter $\frac{\lambda}{1+\lambda}$.
The same holds for $a_3-a_2$, $a_4-a_3$, etc.
By the law of large numbers, $a_n \sim \frac{1+\lambda}{\lambda} n$.
On the other hand, $x_n \sim \mu^{-1} n$.
Therefore, if $\mu< \frac{\lambda}{1+\lambda}$, the probability that $a_k < x_k$ for all $k$ is positive.
This event, in turn, implies success of the algorithm.
\end{proof}

\section{Higher-dimensional arguments}

Let $X^z=(X^z_n)_{n=0,1,2,\dots}$ denote a random walk starting at $z\in\Z^2$ and jumping according to $p(\cdot)$, the superindex is omitted when $z=\o$.

\begin{framed}
\cite{shellef-10}.
$\mu_c \leqslant 1$.
\end{framed}

\begin{proof}
\marginnote{Uses IDLA ideas from Lawler et al. (1992)}
Let $V_n$ be the discrete ball of radius $n$ around $\o$.
Label and order all the particles initially found in $V_n$.
We move the first particle until it is alone in a site, or until it exits the box.
This can be done by subsequently toppling the site that contains the particle until one of these two conditions are met.
If the particle is initially alone, there is nothing to do.
After we are done with the first particle, we move the second particle until it is alone or exits the box.
Then we move the third particle, and so on.

The path performed by each particle in this procedure is distributed as a random walk truncated at the time when it reached an empty site (we do not care about what the particle does once it exits $V_n$).
In order to complete the path, we sample a random walk starting from the point where the particle was stopped.
This second part of the path, which is not actually performed by the particle, is called a \emph{ghost}.
See Figure~\ref{fig:idlash}.

Let $W$ count the \emph{total} number of walks that visit $\o$ before leaving $V_n$.
Let $L$ count the number of walks that visit $\o$ before leaving $V_n$, but which do so \emph{as ghosts}.
We are interested in counting the number of such walks that visit $\o$ \emph{as particles}, which is given by $W-L$.

It may be complicated to determine the set of locations where ghosts start.
But we can use the fact that at most one ghost starts at each site of $V_n$.
Let us start an \emph{artificial ghost} from each site where no actual ghost has been started, as also illustrated in Figure~\ref{fig:idlash}, so that exactly one ghost (artificial or not) starts from each site.
Let $\tilde L \geqslant L$ denote total number of ghosts which visit $\o$ before leaving $V_n$.

\begin{figure}[ht!]
\centering
\includegraphics[width=120mm]{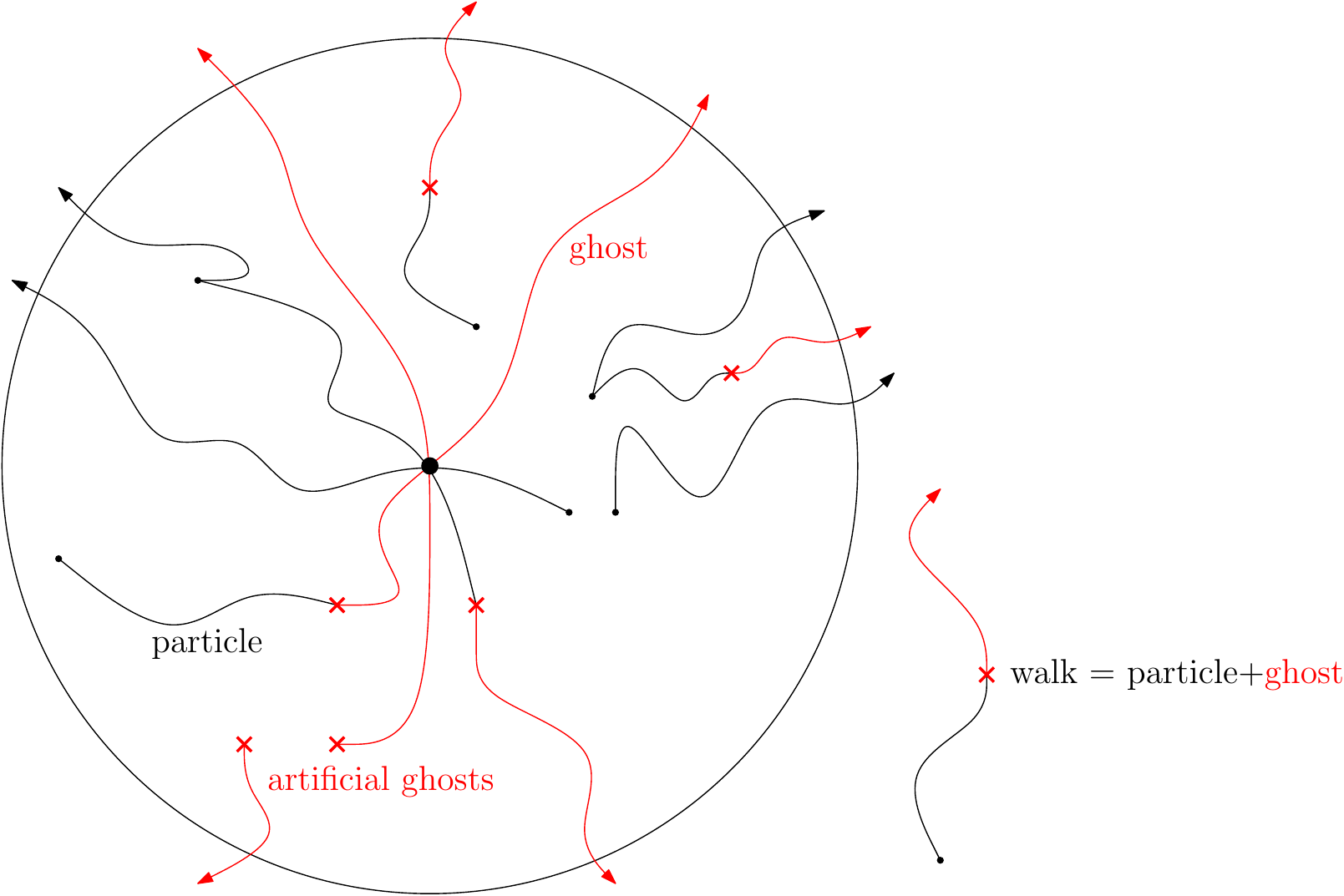}
\caption{Particles which visit $\o$ before or after exiting $V_n$, ghosts started each time a particle find itself alone at a site, and artificial ghosts started from other sites.}
\label{fig:idlash}
\end{figure}

Although $W$ and $\tilde L$ are far from being independent, we claim that each of them has variance bounded by a multiple of its expectation, and thus is concentrated around its mean.
Finally, we also claim that
\[
\E[W] = \mu \E[\tilde L] \to \infty \text{ as } n\to\infty.
\]
Hence, for $\mu>1$, the weak law of large numbers gives $W-L \to \infty$ in probability as $n\to\infty$, proving non-fixation.

To complete the proof, let us prove the above elementary claims.
First, $$\E[\tilde L] = \sum_{x\in V_n} q_x, \quad \text{where} \quad q_x = \Pb[X^x \text{ visits } \o \text{ before exiting } V_n].$$
Also, $\E[W|\eta_0]=\sum_x \eta_0(x) q_x$, so $\E[W] = \E\E[W|\eta_0] = \mu \E[\tilde L]$.
Write $m=\lceil n/3 \rceil$. Estimating the above sum gives
\begin{align*}
\E[\tilde L]
& \geqslant
\sum_{x\in V_{m}} \Pb[X^x \text{ visits } \o \text{ before exiting } V_{m}+x]
\\
& =
\sum_{z\in V_m} \Pb[X^\o \text{ visits } z \text{ before exiting } V_m]
\\
& =
\E[ \text{number of sites visited by } X^\o \text{ before exiting } V_m]
\to \infty
\end{align*}
as $n\to\infty$.
For the variance, first notice that $\tilde L$ is given by the sum of independent indicator functions, so $\V \tilde L \leqslant \E \tilde L$.
Conditioned on $\eta_0$, $W$ is also given by the sum of independent indicator functions, so $\E \V [W|\eta_0] \leqslant \E \E [W|\eta_0] = \E W$ and $\V \E [W|\eta_0] = \sum_x \V[\eta_0(x)] q_x = \mu \mu^{-1} \E[W]$, where $\mu$ denotes first the variance and then the expectation of a Poisson.
Finally, the law of total variance gives $\V[W] \leqslant 2 \E[W]$, finishing the proof.
\end{proof}

\begin{framed}
\cite{taggi-14}
(Biased walks). If the the walk is sufficiently directed, then there exist $\lambda>0$ and $\mu<1$ for which there is no fixation.
\end{framed}

\begin{proof}
Denote the drift by $\v = \sum_z p(z) z \in \R^d \setminus \{\o\}$.
We assume that $\v \cdot \e_1 > 0$ by symmetry, and that $d=2$ for simplicity.
We say that the walk $X^z$ is \emph{good} if $(X^z_n-z) \cdot \e_1 > 0 \text{ for all } n \geqslant 1$.
Let $$\K = \Pb[ (X_n)_n \text{ is good}]>0.$$
Notice that $\K$ tends to $0$ as the walk becomes less biased, and tends to $1$ when the probability of jumping to the left is small.
We will show that the system stays active if
\[
\big[ \mu - \textstyle \frac{\lambda}{1+\lambda}(1-e^{-\mu}) \big] \K > e^{-\mu}.
\]
This condition is meaningful when $\K>\frac{1}{e}$.\footnote
{For Bernoulli initial distribution, the result is a bit stronger. The term $e^{-\mu}$ is replaced by $1-\mu$. In this case, the above condition gives $\mu_c<1$ for any $\lambda<\infty$ and any $\K>0$.}
The result is weaker than in $d=1$ because the argument requires some control on the position where particles exit a large box.

Let $V_n = \{ z=(x,y)\in\Z^2 : x\in[-n,1],y\in[-n^3,n^3] \}$.
Order these sites by writing $V_n=\{z_1,z_2,\dots,z_N\}$, where $N=2n^4$, observing that sites with smaller $x$ must appear first.
We now stabilize $V_n$ as follows.

Pick one particle initially present in $z_1$, and move that particle until it either sleeps or exits $V_n$. If it sleeps in a given site of $V_n$, we leave that particle in that site and do not touch it anymore. In this case we start a \emph{ghost} from that site so that the whole path associated to that particle is distributed as a random walk, part of which is performed by the particle, and part of it by the ghost.
We say that the particle is good if the corresponding walk is good.

Repeat the same procedure for all particles initially present in $z_1$, except for the last particle.
When a single particle is left at $z_1$, its first action may be to sleep, otherwise it is to jump.
If it jumps we treat it as the previous particles, so that we assign a walk to it (part of which may consist of a ghost walk).
If it sleeps, we declare that particle \emph{dead}, and do not even assign it a walk.

After we are finished with all particles initially present in $z_1$, we repeat the same action with particles found on $z_2$, then $z_3$, and so on.

Sites which contain a dead particle are those initially containing at least one particle and such that the last particle sleeps instead of performing its first jump, so the density of dead particles is $(1-e^{-\mu})\frac{\lambda}{1+\lambda}$.
The total number of walks is given by the number of particles minus the number of dead particles, so it has density $$\tilde\mu := \mu - \textstyle \frac{\lambda}{1+\lambda}(1-e^{-\mu}).$$
When the particle is good, its walk is always to the right of its initial position, and the particles initially present on this region have not yet been moved. Therefore, a good particle can only sleep at the set of sites that were initially empty, which has density $e^{-\mu}$.

These observations explain the requirement $\tilde{\mu}\K > e^{-\mu}$.
Let us start a ghost walk at every site that was initially empty, even if no particle was actually stopped there.
This way we obtain an upper bound for the set of ghosts associated to good particles.
It thus follows by the LLN that at least
$( \tilde{\mu} \K - e^{-\mu} ) N - o(N)$ good particles exit $V_n$ before sleeping.
However, we want to estimate the number of particles that exit $V_n$ through site $\o$.

Let $W$ count the number of good walks that exit $V_n$ through $\o$ in this procedure.
Let $L$ count the number of good walks that exit $V_n$ through $\o$, but which do so as ghosts.
We are interested in the number $W-L$ of good particles that exit $V_n$ through $\o$ before sleeping.
In principle it is hard to estimate $L$, so we consider the number $\tilde L \geqslant L$ of ghosts that exit $V_n$ through $\o$, when ghosts are started from each $z\in V_n$ that was initially empty.

Since both $W$ and $\tilde L$ are given by sums of independent variables whose variance is comparable to the expectation (see the proof of~\cite{shellef-10} above), it suffices to show that
\[
{\E[W]}
\geqslant
\tilde{\mu} \K n - o(n)
\qquad
\text{and}
\qquad
{\E[\tilde{L}]}
\leqslant
e^{-\mu} n
\]
in order to have $W-\tilde{L} \to \infty$ as $n\to\infty$, in probability, which in turn implies non-fixation.

We use a simple argument of symmetry.
First, writing $G_k=\{(x,y)\in\Z^2:x=k\}$,
\begin{align*}
\E[\tilde{L}]
&=
\sum_{x=-n}^{-1} \sum_{y=-n^3}^{n^3} e^{-\mu}
\Pb[X^{x,y} \text{ exits } V_n \text{ through } \o]
\\
& \leqslant
e^{-\mu}
\sum_{x=-n}^{-1} \sum_{y\in\Z}
\Pb[X^{x,y} \text{ reaches } G_0 \text{ at } \o]
\\
& =
e^{-\mu}
\sum_{x=-n}^{-1} \sum_{y\in\Z}
\Pb[X^{x,0} \text{ reaches } G_0 \text{ at } (0,-y)]
= e^{-\mu} n
.
\end{align*}
Finally, writing $H_k=\{(x,y)\in\Z^2:y=\pm k\}$,
\begin{align*}
\E[W]
&=
\sum_{x=-n}^{-1} \sum_{y=-n^3}^{n^3} \tilde{\mu}
\Pb[X^{x,y} \text{ is good and exits } V_n \text{ through } \o]
\\
& \geqslant
\tilde{\mu}
\sum_{x=-n}^{-1} \sum_{y=-n^2}^{n^2}
\Pb[X^{x,y} \text{ is good and reaches } G_0 \text{ at } \o \text{, before reaching } y+H_{n^2}]
\\
& =
\tilde{\mu}
\sum_{x=-n}^{-1} \sum_{y=-n^2}^{n^2}
\Pb[X^{x,0} \text{ is good and reaches } G_0 \text{ at } (0,-y) \text{, before reaching } H_{n^2}]
\\
& =
\tilde{\mu}
\sum_{x=-n}^{-1}
\Pb[X^{x,0} \text{ is good and reaches } G_0 \text{ before } H_{n^2}]
\\
& \geqslant
\tilde{\mu}
\sum_{x=-n}^{-1}
\Pb[X^{x,0} \text{ is good and reaches } G_{x+n} \text{ before } H_{n^2}]
\\
& =
\tilde{\mu} n \Pb[X \text{ is good and reaches } G_{n} \text{ before } H_{n^2}].
\end{align*}
Since $\frac{X_k}{k}\to \v$ as $k\to\infty$, the last probability tends to $\K$ as $n\to\infty$.
This proves the desired inequalities, thus finishing the proof.
\end{proof}

\begin{framed}
\cite{shellef-10}.
For $\lambda=\infty$, $\mu_c>0$.
\end{framed}

\begin{proof}
We will in fact show that, if $\mu$ is small, then for almost every initial configuration $\eta_0$, for any realization of $\II$, only finitely many particles can visit the origin.
The idea is to assume that an adversary is trying to bring as many particles to the origin as possible, being able to move the particles to any nearest neighbor rather than following the random jumps.
The proof uses the percolative structure of low-density sets.

Given $\eta_0$, for each finite set $V \subseteq \Z^d$, define its \emph{weight} by $w(V)=\sum_{x\in V} \eta_0(x)$. A finite connected set $V \subseteq \Z^d$ is \emph{internally fillable} if $w(V) \geqslant |V|$.

Consider the decomposition of $\Z^d$ into clusters $\C$ given by
$$\C(x) = \bigcup \big\{ V \subseteq \Z^d : V \text{ is finite and connected, } x \in V, \text{ and }  w(V)\geqslant |V| \text{ or } |V|=1 \big\}.$$
These are called clusters because $\C(x)=\C(y)$ whenever $y\in \C(x)$.

The main observation is that, regardless of $\II$, for a particle starting at $y$ to visit site $x$, it is necessary that some finite connected set $V$ containing both $x$ and $y$ is internally fillable. In other words, it is necessary that $y \in \C(x)$.

Therefore, in order to show fixation, it suffices to prove that, with positive probability, $|\C(x)|<\infty$.
This condition is equivalent to
\[
\sup \Big\{ |V| : \frac{w(V)}{|V|} \geqslant 1 \Big\} < \infty.
\]
It follows from results in \cite{martin-02} about ``greedy lattice animals'' that, for some $C>0$,
\[
\mathrm{a.s.}\lim_{n\to\infty} \max_{V: |V| = n} \frac{w(V)}{|V|} \leqslant C\, \E[\eta_0(\o)^{d+1}].
\]
Since we are assuming Poisson initial distribution, the bound can be made less than $1$ by choosing $\mu$ small enough, which finishes the proof.
\end{proof}

\section{A multi-scale argument}

\begin{framed}
\cite{sidoravicius-teixeira-14}
(Unbiased walks).
For any $\lambda>0$, $\mu_c>0$.
\end{framed}

\begin{proof}
[Overview of the proof]
The proof is too long and technical for these notes.
We give an overview of the general strategy, omitting many delicate points.

The main step is to show that an initial configuration restricted to a very large box stabilizes within a slightly larger box, which high probability.
This is proved by recursion on the scale of the box, so the proof looks very little into the details of the actual ARW dynamics.
In a sense, this kind of approach fits to our intuition that no matter how big a defect is, it will only affect a neighborhood of comparable size.

The box at scale $k$ is a cube of length $L_k$, defined as follows.
Let $\gamma=1/10$, $L_0=10^4$ and
\[
L_{k+1} = \lfloor L_k^\gamma \rfloor ^2 L_k.
\]
Notice that $L_k$ increases as a doubly exponential of $k$.
We also let $R_{k+1} = \lfloor L_k^\gamma \rfloor L_k$ as an intermediate scale between $L_k$ and $L_{k+1}$.
In Figure~\ref{fig:stbox} we see an \emph{inner box} $B_{k+1}'$, an \text{intermediate box} and a \emph{full box} $B_{k+1}$ of level $k+1$.
\begin{figure}[ht]
\small
\centering
\includegraphics[width=55mm]{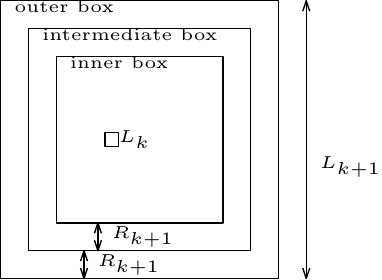}
\caption{Boxes and scales; $L_k \ll R_{k+1} \ll L_{k+1}$.}
\label{fig:stbox}
\end{figure}

Let $p_k$ denote the probability that, starting from a Poisson in $B_k'$, some particle exits $B_k$.
That $p_k \to 0$ fast as $k\to\infty$ follows from the recursion relation
\[
p_{k+1} \leqslant\frac{L_{k+1}^{2d}}{L_k^{2d}} \ {p_k}^2 + e_{k+1},
\]
consisting of a combinatorial term, the probability ${p_k}^2$ that stabilization fails twice at scale $k$, and the probability $e_{k+1}$ that something goes wrong at scale $k+1$.
Indeed, if $p_{k_0}$ is small enough and $e_k \to 0$ fast enough, then the square power above beats the $1+2\gamma$ power in the definition of $L_k$, and $p_k$ vanishes doubly-exponentially fast in $k$.

Let us describe some aspects of this recursion step, depicted in Figure~\ref{fig:renorm}.

\begin{figure}[b!]
\small
\centering
\includegraphics[width=155mm]{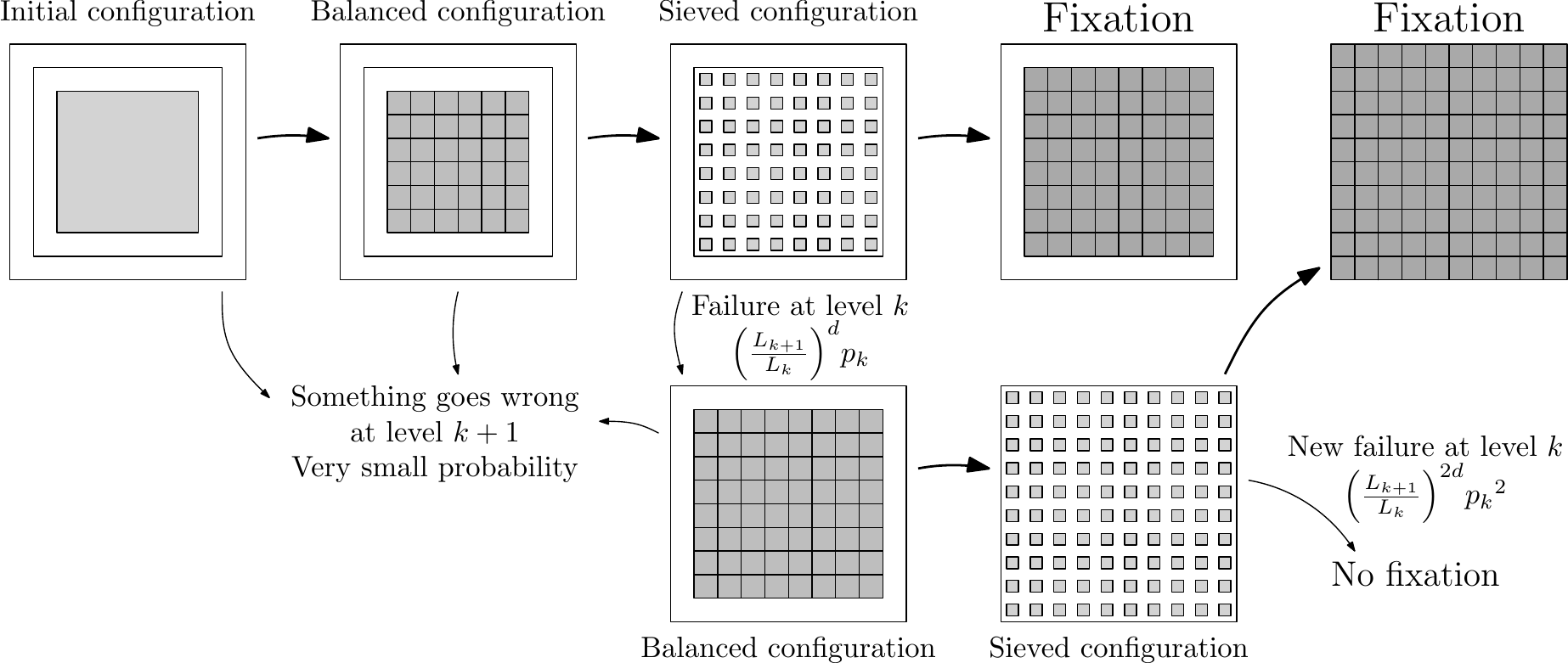}
\caption{Illustrative diagram of events for the recursion relation}
\label{fig:renorm}
\end{figure}

States in light gray have Poisson product distribution with the right density.
They are restricted to the inner box of level $k+1$ for the initial configuration, and the inner boxes of level $k$ for the ``sieved configurations''.
States in dark gray are \emph{stable states}, typically attained by the dynamics.
States in gray are ``balanced configurations''.
Thick arrows represent typical events, while thin arrows represent events of low probability, either $e_{k+1}$ or $p_k$.

Starting fresh. To let the dynamics run on boxes of the previous scale and use recursion, it is important to start with a Poisson product distribution within their inner boxes. This is achieved by a sieving procedure described below.

Worst case scenario.
If the dynamics fail to stabilize all the $\frac{L_{k+1}^d}{L_k^d}$ boxes of level $k$, the configuration inside these boxes is no longer Poisson. In the absence of any useful knowledge about the resulting distribution of particles in this case, we use only the fact that the total number of particles within each box is still a Poisson random variable, and thus cannot be much larger than its mean. A \emph{balanced configuration} is such that the number of particles within each box of level $k$ is appropriately bounded.

Sieving procedure.
Starting from a balanced configurations, we move each particle for quite enough time until their relative position within the box of level $k$ that contains it has mixed. If the particle happens not to be in its inner box, we repeat the procedure until it is. This \emph{reshuffling} with \emph{sieving} results in a state that can be coupled with an i.i.d.\ Poisson configuration with high probability.
This is one of the heaviest statements in~\cite{sidoravicius-teixeira-14} and not simple to prove.
In order for this coupling to be possible, a slight increase in the density is necessary, analogous to the \emph{sprinkling} technique in percolation (this increase should decay just fast enough so that it is summable over $k$).

The chain of events.
By hypothesis we start with a Poisson product measure inside the inner box $B_{k+1}'$.
The first good event is that such a configuration is balanced, that is, each box $B_k$ has a bounded number of particles.
We then let these particles remix by forcing them to jump. During this procedure they can reach the intermediate box, bot not the full box. The second good event is that the resulting configuration is properly sieved.
We now let the evolution run normally within each box $B_k$, and the third good event is that each box stabilizes nicely without letting particles leave.
It may happen that some of these boxes of level $k$ is not stabilized. In this case, the configuration is still balanced with high probability.
We remix them again, now obtaining a sieved configuration in the full box $B_{k+1}$. The system is given a second chance to stabilize, which should happen typically, but may fail again if again some of the level-$k$ boxes does not stabilize as expected.

In the proof there are several aspects to keep under control, and several delicate statements that we omit here.
The above description is not intended to let the reader understand the sketch of proof, but hopefully gives a general flavor of the main argument.
For all the details the reader is referred to the original article.
\end{proof}

\section{Arguments using particle-wise constructions}

The techniques presented so far used the Diaconis-Fulton construction and its properties, as described in~\cite{rolla-sidoravicius-12}.
In the particle-wise construction, the randomness of the jumps is not attached to the sites, but to the particles.
In this section we consider techniques which use the particle-wise construction, or sometimes a combination of both.

\subsection{Preliminaries}

\paragraph{System with labeled particles}
Each existing particle at time $t=0$ is assigned a label $(x,j)$, where $x\in\Z^d$ denotes its starting position and $j=0,\dots,\eta_0(x)-1$ distinguish particles starting at the same site $x$.
Let $Y^{x,j}=(Y^{x,j}_t)_{t \geqslant 0}$ be given by the position of particle $(x,j)$ at each time $t$.
We write $Y^{x,j}_t=\Upsilon \not\in \Z^d$ if $\eta_0(x) \geqslant j$.
Let $\gamma^{x,j}=(\gamma^{x,j}(t))_{t \geqslant 0}$ be given by $\gamma^{x,j}(t) = 1$ if particle $(x,j)$ is active at time $t$ or $\varrho$ if it is passive.
Write $\Y=(Y^{x,j})_{x,j}$ and $\boldsymbol{\gamma}=(\gamma^{x,j})_{x,j}$.
Then the triple $\boldsymbol{\zeta}=(\eta_0,\Y,\boldsymbol{\gamma})$ describes the whole evolution of the system.

Whereas the process $(\eta_t)_{t \geqslant 0}$, given by
$$\eta_t(z)=\sum_x \sum_{j<\eta_0(x)} \gamma^{x,j}(t) \cdot \delta_{Y^{x,j}(t)}(z),$$
only counts the number of particles at a given site at a given time, having each particle labeled gives a lot more information and allows different techniques to be employed.

In a system whose initial configuration $\eta_0$ contains finitely many particles, the evolution described above is always well defined, since it is simply a continuous-time Markov chain on a countable state-space.
Many different constructions will produce $\boldsymbol{\zeta}$ with the correct distribution.

\paragraph{Particle-wise randomness}

Assign to each particle $(x,j)$ a continuous-time walk $X^{x,j}=(X^{x,j}_t)_{t \geqslant 0}$, independently of anything else, as well as a Poisson clock $\PPP^{x,j}\subseteq \R_+$ according to which the particle will try to sleep.
$X^{x,j}$ is the path of the particle parameterized by its \emph{inner time}, which may be slowed down with respect to the system time, depending on the interaction with other particles.
(In the ARW model, when the particle is passive, the inner time halts until it is re-activated by another particle.)
For this reason, $X^{x,j}$ will be called the \emph{putative trajectory} of particle $(x,j)$.
These random elements will be denoted by $\boldsymbol{\xi}=(\eta_0,\X,\boldsymbol{\PPP})$.
Let $\sigma^{x,j}(t)$ denote the inner time of particle $(x,j)$ at instant $t$.
Then $(Y^{x,j}_t)_{t \geqslant 0}$ is given by $Y^{x,j}_t=X^{x,j}_{\sigma^{x,j}(t)}$.

\paragraph{Construction of the infinite system}

When $\eta_0$ contains finitely many particles only, $\boldsymbol{\zeta}$ is determined by $\boldsymbol{\xi}$ in the obvious way, and it works for a.e.\ $\X$ and $\boldsymbol{\PPP}$.
The simultaneous construction on the whole space, is done on the sequences of balls $B(y,n)$, $n\in\N$, centered at each site $y\in\Z^d$.
This family of sequences is countable and translation-invariant.
For each ball, consider the evolution of the system with initial configuration $\eta_0 \cdot \I_{B(y,n)}$.
We say that this construction is \emph{well defined} if (i) for each $x$ and $j$, and for each $T>0$, both $Y^{x,j}_{|_{[0,T]}}$ and $\gamma^{x,j}_{|_{[0,T]}}$ are the same in the systems $(\eta_0 \cdot \I_{B(y,n)},\Y,\boldsymbol{\gamma})$ for all but finitely many $n$, and (ii) the limiting process does not depend on $y$.

The above definition is convenient for two reasons.
First, the family of finite approximations to $\Z^d$ is countable, which makes it possible to prove that the construction is a.s.\ well defined.
Second, this family is translation invariant, so that, whenever the construction is well defined, it is translation covariant.
More precisely, for any translation $\theta$ of $\Z^d$, if $\boldsymbol{\zeta}(\boldsymbol{\xi})$ is well defined, then $\boldsymbol{\zeta}(\theta \boldsymbol{\xi})$ is well defined and equals $\theta \boldsymbol{\zeta}(\boldsymbol{\xi})$, where $\theta \boldsymbol{\xi}=(\theta \eta_0,\theta \X,\theta \boldsymbol{\PPP})$, etc.
In particular, this implies that the system is ergodic, satisfies the mass transport principle, and can be approximated by finite systems regardless of which construction is used.

\paragraph{The Mass-Transport Principle}

Let $m:\Z^d \times \Z^d \to \R_+$ be a translation-invariant random function, that is, $m(x,y;\boldsymbol{\xi})=m(\theta x,\theta y;\theta\boldsymbol{\xi})$ for any translation $\theta$.
The Mass Transport Principle is given by
\[
\E\big[\sum_y m(x,y)\big]=\E\big[\sum_y m(y,x)\big]
.
\]
Informally, the MTP says that the amount of mass transmitted from a vertex $x$ is equal to the amount of mass entering $x$.
It seems like an innocent and perhaps obvious identity, but its strength lies in its versatility, since it holds for \emph{any} such function.
The proof consists simply on re-indexing the sum and using translation invariance. 
See~\cite[Chapter~8]{lyons-peres-} for applications and generalizations to other settings.

\paragraph{The particle-hole model}
Let us introduce a related model that will be useful in the sequel.
Particles perform continuous-time random walks independently of each other.
Sites not containing any particle are called \emph{holes}.
When a particle is alone at some site, it \emph{settles} there forever, filling the corresponding hole.
After the hole has been filled, the site becomes available for other particles to go through.
If a site is occupied by several particles at $t=0^-$, we choose one of them uniformly to fill the hole at $t=0$, and the other particles remain free to move.

As in the ARW with $\lambda=\infty$, once a site has at least one particle, it will always retain one particle.
The differences are (i) sites with $n$ particles are toppled at rate $n$ instead of $n-1$ and (ii) in the system with labeled particles, it is the first particle to arrive at a site that is retained, whereas in the ARW the particles can replace each other.
Nevertheless, both models have the same Diaconis-Fulton representation, and therefore they are equivalent in terms of fixation at $x=\o$.
The Poisson clocks $\boldsymbol{\PPP}$ are not used in either system.

\subsection{Results}

Fixation as defined so far concerns the state of sites, and will be called \emph{site fixation}.
When each labeled particle eventually fixates we call that \emph{particle fixation}.

\begin{framed}
\cite{amir-gurelgurevich-10}
Assume that the particle-wise construction is a.s.\ well defined.
If some particles are not fixating, then sites are not fixating.
Therefore, $\mu_c \leqslant 1$.
\end{framed}

Let us see how $\mu_c \leqslant 1$ follows from the main claim.
First, it follows from the MTP that:
{\par\centering
\emph{The density of particles that fixate cannot be larger than 1.}
\par}
Indeed, let $A(x,j,y)$ denote the event that particle $(x,j)$ fixates at site $y$.
Let $m(x,y)=\sum_j \I_{A(x,j,y)}$.
The density of particles that fixate is given by $\E\sum_y m(x,y)$ and the density of sites eventually occupied by a fixating particle is given by $\E\sum_y m(y,x)$.
These are equal by the MTP.
Since at most one particle can fixate at a given site, $\sum_y m(y,x) \leqslant 1$, proving the estimate.
Now, if $\mu>1$, there is a positive density of particles which do not fixate, which implies that there is no site fixation.
Therefore, $\mu_c \leqslant 1$.

\begin{proof}
Let $A^{v,j}$ denote the event ``$\eta_0(v)>j$ and $(v,j)$ is non-fixating'' and write $A^{v}=A^{v,0}$.
Suppose the probability that some particle is non-fixating is positive. Then there is $j$ such that $\Pb(A^{\o,j})>0$.
By interchangeability of particles, $a:=\Pb(A^\o)=\Pb(\eta_0(\o)>0, (\o,0) \text{ non-fixating}) \geqslant \Pb(\eta_0(\o)>j, (\o,0) \text{ non-fixating}) = \Pb(A^{\o,j})>0$.

As a warm up, notice that by the MTP the number $M_t$ of non-fixating particles present at $\o$ at time $t$ satisfies $\E M_t \geqslant a$, for any $t>0$.
Hence, the $\liminf_{t} \E M_t \geqslant a$.
If we had some control on $\V M_t$, we could conclude that $\liminf_t \Pb(M_t \geqslant 1)>0$, implying non-fixation.
To control variance, we will use a local condition instead, plus extra randomization.

Since the system $\boldsymbol{\zeta}$ is a measurable function of the randomness $\boldsymbol{\xi}$, for any $\epsilon>0$ there is $k\in\N$ such that the event $A^\o$ can be $\epsilon$-approximated by some event $B^\o$ that depends only on $(\eta_0(x),X^{x},\PPP^{x})$ for $\|x\| \leqslant k$.
Let $B^v$ denote the corresponding translation of the event $B^\o$.
When $B^v$ occurs, we say that particle $(v,0)$ is a \emph{candidate}.
It is a \emph{good candidate} if $A^v$ also occurs, otherwise it is a \emph{bad candidate}.

Fix $T>0$.
For $n\in\N$, let us add more randomness to the system by sampling $Z^{v}$ uniformly amongst the first $n$ different sites in the putative trajectory $X^{v,0}$ after time $T$, independently over $v$.
Define $C(v,u)$ as the event ``$B^v$ and $Z^v=u$''.
Let
\[
q(v,u)=\Pb \left( C(v,u) \,\big|\, \boldsymbol{\xi} \right)
\qquad
\text{and}
\qquad
Q(u)=\sum_v q(v,u).
\]
By the mass-transport principle,
\[
\E [Q(v)] = \sum_u \Pb( C(u,v) ) = \sum_u \Pb( C(v,u) ) = \Pb(B^v) =: b>0.
\]
Notice that $q(v,u) \leqslant \frac{1}{n}$.
Notice also that $q(v,u)$ and $q(w,u)$ are independent if $|w-u| > 2k$.
Using these two facts, it can be shown that, as $n\to\infty$,
\(
\V [ Q(v) ] \to 0,
\)
and thus
$Q(v) \to b$ in probability.

Let $N(v)=\sum_u \I_{C(u,v)}$ count the number of candidates for which $Z^u=v$.
Then
\[
\Pb \left( N(v)=0 \,\big|\, \boldsymbol{\xi} \right) =
\prod_u (1-q(u,v)) \leqslant e^{-Q(v)} \to e^{-b} \text{ in probability as }
n \to \infty.
\]
Also, let $\tilde{N}(v) = \sum_u \I_{C(u,v) \setminus A^u}$ count the number of bad candidates for which $Z^u=v$.
Then, using the mass-transport principle,
\[
\E[\tilde{N}(v)]
= \sum_u \Pb[{C(u,v) \setminus A^u}]
= \sum_u \Pb[{C(v,u) \setminus A^v}]
= \Pb[ B^v \setminus A^v ]
\leqslant \epsilon.
\]
Let $D^v$ denote the event ``there exists a good candidate $(u,0)$ such that $Z^u=v$''.
Using the two last estimates we get
\[
\Pb(D^v) \geqslant \Pb(N(v) \geqslant 1) - \Pb(\tilde{N}(v) \geqslant 1 ) \geqslant 1 - e^{-a+\epsilon} -\delta_n - \epsilon,
\]
where $\delta_n \to 0$ as $n\to\infty$.
Choosing $\epsilon$ small and $n$ large, we have $\Pb(D^\o)>\frac{a}{2}$.

On the event $D^\o$, there is some non-fixating particle $(u,0)$ and some inner time $t>T$ such that $X^{u,0}_t = \o$, so vertex $\o$ is visited by an active particle after time $T$.
Letting $T\to\infty$, we get $\Pb(\text{site } \o \text{ not fixating}) \geqslant \frac{a}{2}$, and by the $0$-$1$ law $\Pb(\o \text{ not fixating}) = 1$.
\end{proof}

\begin{framed}
\cite{cabezas-rolla-sidoravicius-14}
For $\lambda=\infty$, $\mu_c \geqslant 1$.
\end{framed}

We start with the following observation about the particle-hole model:
{\par\centering
\emph{The density of holes filled by time $t$ equals the density of particles settled by time $t$.}
\par}
To see that, let $A(x,j,y)$ denote the event that particle $(x,j)$ settles at site $y$ by time $t$, and let $m(x,y)=\sum_j \I_{A(x,j,y)}$.
The density of particles settled by time $t$ is given by $\E\sum_y m(x,y)$ and
the density of holes filled by time $t$ is given by $\E\sum_y m(y,x)$.
These are equal by the MTP, proving the identity.

\begin{proof}
Suppose that $\mu<1$.
Using the above identity,
\[
\Pb[\boldsymbol{o} \mbox{ contains an unfilled hole at time } t] \geq 1-\mu>0.
\]
Since the above event is decreasing in $t$ and the lower bound is uniform with respect to $t$,
\[
\Pb[\boldsymbol{o} \mbox{ is never visited}] > 0.
\]
In particular, the probability that $\o$ is visited finitely many times in the particle-hole model is positive, so by the $0$-$1$ law, $\Pb[m_{\eta_0}(\o)<\infty]=1$, which means that the ARW with $\lambda=\infty$ fixates.
\end{proof}

\begin{framed}
\cite{cabezas-rolla-sidoravicius-14}
No fixation when $\mu=1$.
\end{framed}

The proof is done in full details details in~\cite{cabezas-rolla-sidoravicius-13b} for two-type annihilating systems, that is, with reaction $A+B\to\emptyset$.
Active particles are particles of type $A$, and holes are particles of type $B$.
The reaction happens when a particle settles, thereby filling a hole.

\begin{proof}
[Sketch of the proof]
Suppose that sites fixate in the particle-hole model.
By the first result presented in this section, site fixation implies that every particle eventually settles.
Letting $t\to\infty$ in the previous observation,
\[
\Pb[ \o \text{ is ever visited} ]
=
\mu
.
\]

In the sequel we will show that, under the assumption of site fixation,
\[
\Pb\left[ \o \mbox{ is never visited} \right] > 0.
\]
This implies that $\mu<1$, therefore proving that there cannot be fixation for $\mu=1$.

Assuming site fixation, necessarily, there exists $k\in \N$ such that $$\Pb[\text{the number of particles which ever visit } \o \text{ equals }k]>0.$$
Moreover, there exist $x_1,\dots,x_k\in\Z^d$ such that $\Pb[\mathcal{U}]>0$, where
$$\mathcal{U} = \left[ \text{the particles which ever visit } \o \text{ are initially at the sites }x_1,\dots,x_k \right].$$

Consider two copies $\boldsymbol{\xi}$ and $\tilde{\boldsymbol{\xi}}$ of the system, coupled as follows.
We take $\tilde \X=\X$, and $\tilde\eta_0(x)=\eta_0(x)$ for $x \not \in \{x_1,\dots,x_k\}$.
Finally, for $x \in \{x_1,\dots,x_k\}$, we sample $\tilde\eta_0$ and $\eta_0$ independently.
Now notice that
\begin{multline*}
\Pb\left[ \mathcal{U} \mbox{ occurs for } \tilde{\boldsymbol{\xi}}, \mbox{ and } \eta_0(x_1)=\cdots=\eta_0(x_k)=0 \right]
=
\\
=
\Pb\left[ \mathcal{U} \mbox{ occurs for } \tilde{\boldsymbol{\xi}} \right] \times \Pb\left[ \big. \eta_0(x_1)=\cdots=\eta_0(x_k)=0 \mbox{ for } \boldsymbol{\xi} \right]
>0
.
\end{multline*}
On the above event, no particle ever visits $\o$ in the system $\boldsymbol{\xi}$.
Indeed, on the above event, the initial configuration of $\boldsymbol{\xi}$ is the same as that of $\tilde{\boldsymbol{\xi}}$ except for the deletion of the particles present in $\{x_1,\dots,x_k\}$.
In particular, all the particles which visit the origin in $\tilde{\boldsymbol{\xi}}$ are deleted in $\boldsymbol{\xi}$.
Recalling that $\boldsymbol{\xi}$ and $\tilde{\boldsymbol{\xi}}$ share the same putative trajectories, by following how the effect of deleting such particles propagates in the system evolution, we see that in the system $\boldsymbol{\xi}$ no particles can ever visit $\o$.
\end{proof}

\bibliographystyle{bib/leo}
\addcontentsline{toc}{section}{References}
\bibliography{bib/leo}

\end{document}

%% file: macros.tex
\usepackage{amscd}
\usepackage{amsfonts}
\usepackage{amsmath}
\usepackage{amssymb}
\usepackage{amstext}
\usepackage{amsthm}
\usepackage{color}
\usepackage[nodate]{datetime}
\usepackage{dsfont}
\usepackage{framed}
\usepackage{graphicx}
\usepackage{hyperref}
\usepackage[latin1]{inputenc}
\usepackage{mathrsfs}
\usepackage{mathtools}
\usepackage{srcltx}
\usepackage{stmaryrd}
\usepackage[normalem]{ulem}
\usepackage{verbatim}
\usepackage{vruler}
\usepackage{wasysym}

\usepackage[T1]{fontenc}
\usepackage{lmodern}

\usepackage[lmargin=25mm,rmargin=25mm,top=25mm,bottom=30mm,paperwidth=210mm,paperheight=11in]{geometry}

\newcommand{\mydate}{{\today}}

\date{\mydate}

\newcommand{\marginnote}[1]{\marginpar{\raggedright\tiny #1}}
\renewcommand{\marginnote}[1]{}

\newtheorem*{theorem*}{Theorem}

\newtheorem*{observation*}{Observation}
\theoremstyle{definition}

\newtheorem*{remark*}{Remark}

\newcommand{\E}{\mathbb{E}}
\newcommand{\N}{\mathbb{N}}
\newcommand{\Pb}{\mathbb{P}}
\newcommand{\PPP}{\mathcal{P}}
\newcommand{\R}{\mathbb{R}}
\newcommand{\V}{\mathbb{V}}
\newcommand{\Z}{\mathbb{Z}}

\newcommand{\C}{\mathcal{C}}

\newcommand{\K}{\mathcal{K}}

\newcommand{\X}{\mathbf{X}}
\newcommand{\Y}{\mathbf{Y}}

\renewcommand{\geq}{\geqslant}

\newcommand{\h}{h}
\newcommand{\e}{{\boldsymbol{e}}}
\renewcommand{\o}{{\boldsymbol{o}}}
\renewcommand{\o}{{\boldsymbol{0}}}
\renewcommand{\v}{{\boldsymbol{v}}}
\newcommand{\I}{\mathds{1}}
\newcommand{\II}{\mathcal{I}}